\documentclass{amsart}

\usepackage{amssymb}
\usepackage{graphicx}

\newtheorem{theorem}{Theorem}[section]

\newtheorem{proposition}[theorem]{Proposition}
\newtheorem{corollary}[theorem]{Corollary}

\newtheorem{_algorithm}[theorem]{Algorithm}

\newtheorem{_definition}[theorem]{Definition}
\newenvironment{definition}{\begin{_definition}\rm}{\end{_definition}}

\newtheorem{_remark}[theorem]{\it Remark}
\newenvironment{remark}{\begin{_remark}\rm}{\end{_remark}}

\newtheorem{_example}[theorem]{Example}

\newtheorem{_assumption}[theorem]{Assumption}

\newtheorem{_construction}[theorem]{Construction}

\newtheorem{_claim}[theorem]{Claim}
\newenvironment{claim}{\begin{_claim}\rm}{\end{_claim}}

\newtheorem{_conjecture}[theorem]{Conjecture}
\newenvironment{conjecture}{\begin{_conjecture}\rm}{\end{_conjecture}}

\numberwithin{equation}{section}
\numberwithin{table}{section}
\numberwithin{figure}{section}
\renewcommand{\qed}{\hfill {$\Box$}}

\newcommand{\A}{\mathord{\mathbb A}}

\newcommand{\F}{\mathord{\mathbb F}}
\renewcommand{\P}{\mathord{\mathbb  P}}
\newcommand{\Q}{\mathord{\mathbb  Q}}

\newcommand{\Z}{\mathord{\mathbb Z}}

\newcommand{\UUU}{\mathord{\mathcal U}}

\newcommand{\mapdownsurj}{
\hbox{$\bigm\downarrow$}
\llap{\hbox{\raise 2pt\hbox{$\bigm\downarrow$}}}%
\vstrechmapdown
}

\newcommand{\isom}{\smash{\mathop{\;\to\;}\limits\sp{\sim\,}}}

\newcommand{\set}[2]{\{\; {#1} \; \mid \; {#2} \;  \}}
\newcommand{\shortset}[2]{\{ {#1} \,|\, {#2}   \}}

\newcommand{\sm}{\setminus}
\newcommand{\st}{\subset}

\newcommand{\wt}{\widetilde}

\newcommand{\sprime}{\sp\prime}

\newcommand{\sperp}{\sp{\perp}}

\newcommand{\dual}{\sp{\vee}}

\newcommand{\semidirectproduct}{\rtimes}

\newcommand{\inv}{\sp{-1}}

\newcommand{\NS}{\mathord{\rm NS}}

\newcommand{\Aut}{\operatorname{\rm Aut}\nolimits}
\newcommand{\pr}{\operatorname{\rm pr}\nolimits}
\newcommand{\Sing}{\operatorname{\rm Sing}\nolimits}

\newcommand{\Pt}{\P^2}

\newcommand{\der}{\partial}

\newcommand{\Der}[2]{\frac{\der #1}{\der #2}}

\newcommand{\ang}[1]{\langle #1\rangle}

\newcommand{\rmand}{\textrm{and}}

\newcommand{\quand}{\quad\rmand\quad}

\begin{document}

\title[Supersingular $K3$ surfaces in characteristic $5$]{
Unirationality of certain
supersingular $K3$ surfaces in characteristic $5$}

\author{Duc Tai Pho}
\address{
Department of Mathematics,
Vietnam National University,
334 Nguyen Trai street, Hanoi, VIETNAM
}
\email{phoductai@yahoo.com, taipd@vnu.edu.vn}

\author{Ichiro Shimada}
\address{
Department of Mathematics,
Faculty of Science,
Hokkaido University,
Sapporo 060-0810,
JAPAN
}
\email{shimada@math.sci.hokudai.ac.jp
}


\begin{abstract}
We show that every
supersingular $K3$ surface in characteristic $5$ with Artin invariant $\le 3$
is unirational.
\end{abstract}

\maketitle

\section{Introduction}
We work over an algebraically closed field $k$.
\par
\medskip
A $K3$ surface $X$ is called \emph{supersingular}
(in the sense of Shioda~\cite{MR0374149})
if the Picard  number of $X$ is equal to the second Betti number $22$.
Supersingular $K3$ surfaces exist only when  the characteristic
of $k$ is positive.
Artin~\cite{MR0371899} showed that,
if $X$ is a supersingular $K3$ surface in characteristic $p>0$,
then the discriminant of the N\'eron-Severi lattice $\NS (X)$ of $X$ is
written as  $-p^{2\sigma(X)}$,
where $\sigma(X)$ is a positive integer $\le 10$.
(See also Illusie~\cite[Section~7.2]{MR565469}.)
This integer $\sigma(X)$ is called the \emph{Artin invariant of $X$}.
\par
\medskip
A surface $S$ is called \emph{unirational} if
the function field $k(S)$ of $S$ is contained
in a purely transcendental extension field of $k$,
or equivalently, if there exists a dominant rational map from a projective plane
$\Pt$  to $S$.
Shioda~\cite{MR0374149} proved that,
if a smooth projective surface $S$ is unirational,
then the Picard number of $S$ is equal to the second Betti number of $S$.
Artin and Shioda conjectured
that the converse is true for $K3$ surfaces
(see, for example, Shioda~\cite{MR0429911}):
\begin{conjecture}
Every supersingular $K3$ surface is unirational.
\end{conjecture}
In this paper, we consider this conjecture for
supersingular $K3$ surfaces in characteristic $5$.
\par
\medskip
From now on, we assume that the characteristic of $k$ is $5$.
Let $k[x]_6$ be the space of polynomials in $x$ of degree $6$,
and let $\UUU\subset k[x]_6$ be the space of  $f(x)\in k[x]_6$ 
such that the quintic equation $f\sp\prime (x)=0$ has no
multiple roots.
It is obvious that $\UUU$ is
a Zariski open dense subset  of $k[x]_6$.
For $f\in \UUU$, we denote by  $C_f\subset \P^2$ 
 the projective plane curve of degree $6$
 whose affine part
is  defined by 
$$
y^5-f(x)=0.
$$
Let $Y_f\to\P^2$ be the double covering of $\P^2$ 
whose branch locus is equal to $C_f$, 
and
let $X_f\to Y_f$ be the minimal resolution of $Y_f$.
\begin{theorem}\label{thm:main}
If $f$ is a polynomial in $\UUU$, then $X_f$ is a supersingular $K3$ surface
with $\sigma (X_f)\le 3$.
Conversely, if $X$ is a supersingular $K3$ surface
with $\sigma (X)\le 3$, then there exists $f\in \UUU$ such that
$X$ is isomorphic to $X_f$.
\end{theorem}
The affine part of $Y_f$ is defined by 
$w^2=y^5-f(x)$.
Hence 
the function field $k(X_f)$ is equal to $k(w, x,y)$,
and it is contained in the purely transcendental extension field
$k(w^{1/5}, x^{1/5})$ of $k$.
Therefore we obtain the following corollary:
\begin{corollary}
Every supersingular $K3$ surface in characteristic $5$
with Artin invariant $\le 3$ is unirational.
\end{corollary}
The unirationality of a supersingular $K3$ surface $X$ in characteristic $p>0$
with Artin invariant $\sigma$ has been proved in the following cases:
(i) $p=2$, (ii)  $p=3$ and $\sigma\le 6$,
and (iii) $p$ is odd and $\sigma\le 2$.
In the cases (i) and (ii), the unirationality 
was proved by Rudakov~and~Shafarevich~\cite{MR508830},~\cite{MR633161}
by showing that there exists a structure of the quasi-elliptic fibration on $X$.
The case (iii) follows from the result of Ogus~\cite{MR563467},\cite{MR717616}
that a supersingular $K3$ surface in odd characteristic with Artin invariant $\le 2$
is a Kummer surface associated with a supersingular abelian surface,
and the result of Shioda~\cite{MR0572983} that such a Kummer surface is unirational.
The unirationality of $X$ in the case $(p, \sigma)=(5, 3)$
proved in this paper seems to be  new.
\par
\medskip
In~\cite{MR2059747},
we have shown that a supersingular $K3$ surface in characteristic $2$
is birational to a normal $K3$ surface with $21 A_1$-singularities,
and that such a normal $K3$ surface is a purely inseparable double cover of $\P^2$.
In~\cite{tencusps},
we have proved that a supersingular $K3$ surface in characteristic $3$
with Artin invariant $\le 6$ is birational to a normal $K3$ surface with $10 A_2$-singularities,
and it is also birational to a purely inseparable triple cover of $\P^1\times\P^1$.
These yield an alternative proof
to the results of Rudakov~and~Shafarevich~\cite{MR508830},~\cite{MR633161}
in the cases (i) and (ii) above.
\par
\medskip
In this paper,
we show that 
a supersingular $K3$ surface in characteristic $5$
with Artin invariant $\le 3$ is birational to 
a normal $K3$ surface with $5 A_4$-singularities
that is a double cover of $\P^2$, 
and then prove that such a normal $K3$ surface is 
isomorphic to $Y_f$ for some $f\in \UUU$.
The first step follows from the structure theorem
of the N\'eron-Severi lattices of supersingular $K3$ surfaces
due to Rudakov~and~Shafarevich~\cite{MR633161}.
For the second step,
we investigate
projective plane curves of degree $6$
with $5A_4$-singularities in Section~\ref{sec:curve}.
\section{Projective plane curves with $5A_4$-singularities}\label{sec:curve}
\begin{definition}
A germ of a curve singularity in characteristic $\ne 2$
is called an $A_n$-singularity if it is formally isomorphic to
$$
y^2-x^{n+1}=0,
$$
(see Artin~\cite{MR0450263},  and Greuel~and~Kr\"oning~\cite{MR1033443}.)
\end{definition}
We assume that the base field $k$ is of characteristic $5$  until  the end of the paper.
\begin{proposition}\label{prop:C}
Let $C\subset\Pt$ be a reduced projective plane curve of degree $6$.
Then the following conditions are equivalent to each other.
\begin{itemize}
\item[(i)] The singular locus of $C$ consists of five $A_4$-singular points.
\item[(ii)] There exists $f\in \UUU$ such that $C=C_f$.
\end{itemize}
\end{proposition}
For the proof, we need the following result due to Wall~\cite{MR1317485},
which holds in any characteristic.
Let $D\subset\Pt$ be an integral  plane curve of degree $d>1$,
and let $I_D\subset \Pt\times(\Pt)\dual$ be the closure of the locus
of all $(x, l)\in \Pt\times(\Pt)\dual$ such that
$x$ is a smooth point of $D$ and $l$ is the tangent line to $D$ at $x$.
Let $D\dual\subset(\Pt)\dual$ be the image of the second projection
$$
\pi_D: I_D\to (\Pt)\dual.
$$
We equip $D\dual$ with the reduced structure,
and call it the \emph{dual curve of $D$}.
Note that the first projection $I_D\to D$ is birational.
Therefore, by the projection $\pi_D$,
we can regard the function field $k(D)$ as an extension field of
the function field $k(D\dual)$.
The corresponding rational map 
from $D$ to $D\dual$ is called the \emph{Gauss map}.
We put
$$
\deg \pi_D:=[k(D): k(D\dual)].
$$
We choose general homogeneous coordinates $[w_0:w_1:w_2]$ of $\Pt$,
and let $F(w_0, w_1, w_2)=0$ be the defining equation of $D$.
We denote by  $D_Q\subset \Pt$  the curve defined by
$$
\Der{F}{w_2}=0,
$$
which  is called the \emph{polar curve of $D$} with respect to $Q=[0:0:1]$.
\begin{proposition}[Wall~\cite{MR1317485}]\label{prop:Wall}
For a singular point $s$ of $D$,
we denote by $(D.D_Q)_s$ the local intersection multiplicity of $D$ and 
$D_Q$ at $s$. Then we have
$$
\deg \pi_D \cdot \deg D\dual=d(d-1)-\sum_{s\in \Sing (D)} (D.D_Q)_s.
$$
\end{proposition}
\begin{remark}\label{rem:A4}
If $s\in D$ is an $A_n$-singular point, then 
the polar curve $D_Q$ is smooth at $s$ and 
the local intersection multiplicity $(D.D_Q)_s$ is $n+1$.
\end{remark}
\begin{proof}[Proof of Proposition~\ref{prop:C}]
Suppose that $C$ has $5 A_4$-singular points as its only singularities.
Since an $A_4$-singular point is unibranched,
$C$ is irreducible.
By Proposition~\ref{prop:Wall} and Remark~\ref{rem:A4},
we have 
$$
\deg \pi_C\cdot  \deg C\dual =5.
$$
Suppose that $(\deg \pi_C ,\deg C\dual) =(1, 5)$.
Let $\nu: \wt{C}\to C$ be the normalization of $C$.
Since $\deg \pi_C =1$, we can consider $\wt{C}$
as  a normalization of $C\dual$.
We denote by
$$
\nu\dual : \wt{C}\to C\dual
$$
the  morphism of normalization.
Let $s$ be a singular point of $C$,
and let  $\wt{s}\in \wt{C}$ be the point of $\wt{C}$
that is mapped to $s$ by $\nu$.
We can choose affine coordinates $(x, y)$ of $\Pt$
with the origin $s$ and a formal parameter $t$ of $\wt{C}$ at $\wt{s}\,$
such that $\nu$ is given by
$$
t\mapsto (x, y)=(\,t^2, \; t^5+\,c_6\, t^6+ \,c_7\, t^7+\cdots).
$$
Let $(u, v)$ be the affine coordinates of $(\Pt)\dual$
such that the point $(u, v)\in (\Pt)\dual$
corresponds to  the line of $\Pt$ defined by $y=ux+v$.
Then $\nu\dual$  is given at $\wt{s}\,$ by
$$
t\mapsto (u, v)=(\,3\,c_6\, t^4+\cdots, \;  t^5+\cdots).
$$
(See, for example, Namba~\cite[p.~78]{MR0768929}.)
Therefore $\nu\dual(\wt{s}\,)$ is a singular point of $C\dual$
with multiplicity $\ge 4$.
We choose distinct two points $s_1, s_2\in \Sing (C)$.
There exists a line of $(\Pt)\dual$ that passes through
both of $\nu\dual (\wt{s_1})\in C\dual$ and $\nu\dual(\wt{s_2})\in C\dual$.
This contradicts Bezout's theorem,
because $\deg C\dual =5< 4+4$.
Therefore we have
$(\deg \pi_C ,\deg C\dual) =(5,1)$.
Then there exists a point $P\in \Pt$ such that
we have
\begin{equation}\label{eq:equiv1}
l\in C\dual\;\;\Longleftrightarrow\;\;P\in l.
\end{equation}
We choose homogeneous coordinates $[w_0:w_1:w_2]$ of $\P^2$
in such a way that  $P=[0:1:0]$.
Let $L_\infty$ be the line $w_2=0$, and
let $(x, y)$ be the affine coordinates on $\A^2:=\P^2\setminus L_\infty$
given by $x:=w_0/w_2$ and $ y:=w_1/w_2$.
Suppose that $C$ is defined by $h(x, y)=0$ in $\A^2$.
From~\eqref{eq:equiv1}, we have%
\begin{equation}\label{eq:equiv2}
h(a, b)=0\;\;\Longrightarrow\;\;\frac{\der h}{\der y} (a, b)=0.
\end{equation}
Let $U_C\subset \A^1$ be the image of the projection $(C\sm \Sing(C))\cap \A^2\to \A^1$
given by $(a, b)\mapsto a$.
Note that $U_C$ is Zariski dense in $\A^1$.
Let $(a_0, b_0)$ be a smooth point of $C\cap\A^2$.
By~\eqref{eq:equiv2}, we have
$$
\frac{\der h}{\der x} (a_0, b_0)\ne 0.
$$
Hence 
there exists a formal power series $\gamma (\eta)\in k[[\eta]]$ such that
$C$ is defined by $x-a_0=\gamma (y-b_0)$
locally around  $(a_0, b_0)$.
By~\eqref{eq:equiv2} again,
$\gamma\sprime (\eta)$ is constantly equal to $0$,
and hence there exists a formal power series $\beta(\eta)\in k[[\eta]]$
such that $\gamma (\eta)=\beta(\eta)^5$.
Therefore the local intersection multiplicity of the line $x-a_0=0$ and $C$ at $(a_0,b_0)$
is $\ge 5$.
Thus we obtain the following:
\begin{equation}\label{eq:equiv3}
\parbox{7cm}{If $a\in U_C$, then the equation $h(a, y)=0$ in $y$
has a root of multiplicity $\ge 5$.}
\end{equation}
We put
$$
h(x, y)\;=\;c\,y^6 \,+\, g_1(x)\, y^5\,+\,\cdots\, + \,g_5(x)\,y\,+\,g_6(x),
$$
where $c$ is a constant, and $g_{\nu} (x)\in k[x]$ is a polynomial of degree $\le \nu$.
Suppose that $c\ne 0$.
We can assume $c=1$.
By~\eqref{eq:equiv3}, we have
$g_2(a)=g_3(a)=g_4(a)=0$ and $g_1(a)g_5(a)=g_6(a)$ for any $a\in U_C$.
Since $U_C$ is Zariski dense in $\A^1$,
we have $g_2=g_3=g_4=0$ and $g_1 g_5=g_6$.
Then we have $h(x,y)=(y^5+g_5(x))(y+g_1(x))$,
which contradicts  the irreducibility of $C$.
Thus  $c=0$ is proved.
Then, by~\eqref{eq:equiv3}, we have
$g_1\ne 0$ and $g_2=g_3=g_4=g_5=0$.
We put $g_1=Ax+B$,
and define a new homogeneous coordinate system $[z_0:z_1:z_2]$
of $\Pt$ by
$$
\begin{cases}
(z_0,z_1,z_2):=(w_0, w_1, Aw_0+Bw_2) & \textrm{if $B\ne 0$;} \\
(z_0,z_1,z_2):=(w_2, w_1, Aw_0) & \textrm{if $B=0$.}
\end{cases}
$$
Then $C$ is defined by a homogeneous equation of the form
$$
z_2 z_1^5-F(z_0, z_2)=0,
$$
where $F(z_0, z_2)$ is a homogeneous polynomial of degree $6$.
We put $L\sprime_\infty:=\{z_2=0\}$.
Defining  the affine coordinates $(x, y)$ 
on $\P^2\sm L\sprime_\infty$ by  $(x, y):=(z_0/z_2, z_1/z_2)$,
we see that the affine part of $C$ is defined by $y^5-f(x)$
for some polynomial $f (x)$ of degree $\le 6$.
If $\deg f<6$, then $L\sprime_\infty$ would be an irreducible component of $C$
because $\deg C=6$.
Therefore we have $\deg f=6$.
Then $C\cap L\sprime _\infty$ consists of a single point $[0:1:0]$,
and  $C$ is smooth at $[0:1:0]$.
Therefore we have
$$
\Sing (C)=\set{(\alpha, f(\alpha)^{1/5})}{f\sprime (\alpha)=0}.
$$
Since $C$ has  five singular points,
we have $f\in\UUU$.
\par
\medskip
Conversely, suppose that $f\in \UUU$.
We show that $\Sing (C_f)$ consists of   $5A_4$-singular points.
Let $L_\infty\st \Pt$ be the line at infinity.
It is easy to check that
$C_f\cap L_\infty $ consists of a single point $[0:1:0]$,
and $C_f$ is smooth at this point.
Therefore we have 
$\Sing (C_f)=\shortset{(\alpha, f(\alpha)^{1/5})}{f\sprime (\alpha)=0}$.
In particular,  $C_f$ has exactly five singular points.
Let $(\alpha, \beta)$ be a singular point of $C_f$.
Since $\alpha$ is a simple root of the quintic equation $f\sprime(x)=0$, 
there exists a polynomial $g(x)$ with $g(\alpha)\ne 0$ such that 
$$
f(x)=f(\alpha) + (x-\alpha)^2 g(x).
$$
Because $\beta^5=f(\alpha)$,
the defining equation of $C$ is written as
$$
(y-\beta)^5-(x-\alpha)^2 g(x)=0.
$$
Therefore $(\alpha, \beta)$ is an $A_4$-singular point of $C_f$.
\end{proof}
\section{Proof of Theorem~\ref{thm:main}}
First we show that, if $f\in \UUU$, then $X_f$ is a supersingular $K3$ surface
with Artin invariant $\le 3$.
Since the sextic double plane $Y_f$ has only rational double points as its singularities by
Proposition~\ref{prop:C},
its minimal resolution $X_f$ is a $K3$ surface by the results of Artin~\cite{MR0146182},~\cite{MR0199191}.
Let $\Sigma_f$ be  the sublattice of the N\'eron-Severi lattice $\NS(X_f)$
of $X_f$ that is generated by the classes of the $(-2)$-curves
 contracted by $X_f\to Y_f$.
Then $\Sigma_f$ is isomorphic to the negative-definite root lattice
of type $5A_4$ by Proposition~\ref{prop:C}.
In particular, $\Sigma_f$ is of rank $20$,  and its discriminant is $5^5$.
Let $H_f\subset X_f$ be the pull-back of a line of $\Pt$,
and put
$$
h_f:=[H_f]\in \NS(X_f).
$$
Since the line at infinity $L_\infty\subset \Pt$ intersects $C_f$ at a single point $[0:1:0]$
with multiplicity $6$,
and $[0:1:0]$ is a smooth point of $C_f$,
the pull-back of $L_\infty$ to $X_f$ is a union of two smooth
rational curves that intersect each other at a single point with multiplicity $3$.
Let $L_f$ be one of the two rational curves,
and  put
$$
l_f:=[L_f]\in \NS(X_f).
$$
Then $h_f$ and $l_f$ generate a lattice $\ang{h_f,l_f}$ of rank $2$ in $\NS (X_f)$
whose  intersection matrix is equal to 
$$
\left(
\begin{array}{cc}
2 & 1 \\
1 & -2
\end{array}
\right).
$$
In particular, the discriminant of $\ang{h_f,l_f}$ is $-5$.
Note that $\Sigma_f$ and $\ang{h_f,l_f}$ are orthogonal in $\NS(X_f)$.
Therefore $\NS(X_f)$ contains a sublattice $\Sigma_f\oplus \ang{h_f,l_f}$
of rank $22$ and discriminant $-5^6$.
Thus $X_f$ is supersingular, and $\sigma(X_f)\le 3$.
\par
\medskip
In order to prove the second assertion of Theorem~\ref{thm:main},
we define an even lattice $S_0$ of rank $22$ with signature $(1, 21)$ 
and discriminant $-5^6$ by
$$
S_0:=\Sigma^{-}_{5A_4}\oplus \ang{h,l},
$$
where $\Sigma^{-}_{5A_4}$ is the negative-definite root lattice
of type $5A_4$,
and $\ang{h,l}$ is the lattice of rank $2$ generated by
the vectors $h$ and $l$ satisfying
$$
h^2=2,\quad l^2=-2,\quad hl=1.
$$
\begin{remark}
This lattice $\ang{h,l}$ is the unique even indefinite lattice of rank $2$
with discriminant $-5$.
See Edwards~\cite{MR1416327}, or Conway and Sloane~\cite[Table 15.2a]{MR1662447}.
\end{remark}
\begin{claim}\label{claim:overlattices}
For $\sigma=1, 2, 3$,
there exists an even overlattice $S^{(\sigma)}$ of $S_0$
with the following properties:
\begin{itemize}
\item[(i)] the discriminant of $S^{(\sigma)}$ is $-5^{2\sigma}$,
\item[(ii)] the Dynkin type of the root system
$\shortset{r\in S^{(\sigma)}}{rh=0, r^2=-2}$
is $5A_4$, 
\item[(iii)]
the set
$\shortset{e\in S^{(\sigma)}}{eh=1, e^2=0}$
is empty.
\end{itemize}
\end{claim}
Here we prove that $S^{(3)}=S_0$ satisfies (ii) and (iii).
Let $v=s+xh+yl$ be a vector of $S^{(3)}=S_0$,
where $s\in \Sigma^{-}_{5A_4}$ and $x, y\in \Z$.
If $vh=0$ and $v^2=-2$,
then we have $2x+y=0$ and $s^2-10x^2=-2$.
Since $s^2\le 0$,
we have $x=y=0$ and hence $v$ is a root in $\Sigma^{-}_{5A_4}$.
Therefore $S^{(3)}=S_0$ satisfies (ii).
If $vh=1$ and $v^2=0$,
then we have $2x+y=1$ and $s^2-10x^2+10x-2=0$.
Since $s^2\le 0$,
there is not   such an integer $x$.
Hence $S^{(3)}=S_0$ satisfies (iii).
Thus Claim~\ref{claim:overlattices} for $\sigma=3$ has been  proved.
For the cases $\sigma=2$ and $\sigma=1$, see 
Proposition~\ref{prop:Hs} in
the next section.
%
%
\par
\medskip
Let $X$ be a supersingular $K3$ surface with $\sigma=\sigma (X)\le 3$.
By the results of Rudakov and Shafarevich~\cite{MR633161},
the isomorphism class of the lattice $\NS(X)$ is characterized by the following properties;
\begin{itemize}
\item[(a)] even and signature $(1, 21)$, and
\item[(b)]  the discriminant group is isomorphic to $\F_5^{\oplus 2\sigma}$.
\end{itemize}
Since the discriminant group of $S^{(\sigma)}$
is a quotient group of a subgroup of 
the discriminant group $\F_5^{\oplus 6}$ of $S_0$,
the lattice $S^{(\sigma)}$ has also these properties.
Therefore 
there exists an isomorphism
$$
\phi\;:\; S^{(\sigma)} \;\isom\;  \NS(X).
$$
By~\cite[Proposition 3 in Section 3]{MR633161},
we can assume that $\phi (h)$ is the class $[H]$ of a nef divisor $H$.
Note that $H^2=h^2=2$.
If the complete linear system $|H|$ had  a fixed component, 
then,
by Nikulin~\cite[Proposition 0.1]{MR1260944},
 there would be  an elliptic pencil $|E|$ and a $(-2)$-curve $\Gamma$
such that $|H|=2|E|+\Gamma$ and $ E \Gamma=1$,
and  the vector $e\in S^{(\sigma)}$
that is mapped to $[E]$ by $\phi$
 would satisfy
$eh=1$ and $e^2=0$.
Therefore  
 the property (iii) of $S^{(\sigma)}$
implies that  the  linear system $|H|$ has no fixed components
(see also Urabe~\cite[Proposition 1.7]{MR1101859}.)
Then, by Saint-Donat~\cite[Corollary 3.2]{MR0364263}, $|H|$ is base point free.
Hence  we have a morphism
$\Phi_{|H|}: X\to \Pt$ induced by $|H|$.
Let
$$
X\;\to\; Y_H\;\to\; \Pt
$$
be the Stein factorization of $\Phi_{|H|}$.
Then $Y_H\to \Pt$ is a finite double covering
branched along a curve
$C_H\subset \Pt$ of degree $6$.
By the property (ii) of $S^{(\sigma)}$,
we see that $\Sing (Y_H)$ consists of 
 $5A_4$-singular points,
and hence  $\Sing (C_H)$
also consists of  $5A_4$-singular points.
By Proposition~\ref{prop:C},
there exists an element $f\in \UUU$ such that $C_H$ is isomorphic to $C_f$.
Then $X$ is isomorphic to $X_f$.
\qed
\begin{remark}
In~\cite{Milnor20},
it is proved that a normal $K3$ surface with $5A_4$-singular points
exists only in characteristic $5$.
\end{remark}
\section{Classification of overlattices}
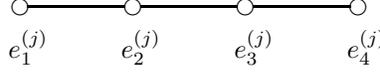
\begin{figure}
\begin{center}
\setlength{\unitlength}{.75mm}
\begin{picture}(60, 17)(0,0)
\put(0,10){\circle{2.8}}
\put(20,10){\circle{2.8}}
\put(40,10){\circle{2.8}}
\put(60,10){\circle{2.8}}
\put(-2,1){$e_1^{(j)}$}
\put(18,1){$e_2^{(j)}$}
\put(38,1){$e_3^{(j)}$}
\put(58,1){$e_4^{(j)}$}
\put(1.4,10){\line(1,0){17.2}} 
\put(21.4,10){\line(1,0){17.2}} 
\put(41.4,10){\line(1,0){17.2}} 
\end{picture}
\end{center}
\caption{The Dynkin diagram of type $A_4$}\label{fig:dynkinA4}
\end{figure}

Let $F\subset S_0$ be a fundamental system of roots
of $\Sigma^{-}_{5A_4}\subset S_0$
(see Ebeling~\cite{MR1938666} for the definition 
and properties of a fundamental system of roots.)
Then $F$ consists of $4\times 5$ vectors 
$$
e_i^{(j)}\quad(i=1, \dots, 4, \;j=1, \dots, 5)
$$
such that
$$
e^{(j)}_i e^{(j\sprime)}_{i\sprime}=
\begin{cases}
0 & \textrm{if $j\ne j\sprime$ or $|i-i\sprime|>1$,} \cr
1 & \textrm{if $j= j\sprime$ and $|i-i\sprime|=1$,}\; \cr
-2 & \textrm{if $j= j\sprime$ and $i=i\sprime$,}
\end{cases}
$$
(see Figure~\ref{fig:dynkinA4}.)
We put
$$
\Aut (F, h):=\set{g\in O(S_0)}{g(F)=F,  g(h)=h},
$$
where $O(S_0)$ is the orthogonal group of the lattice $S_0$.
Then $\Aut (F, h)$ is isomorphic to the automorphism group 
of the Dynkin diagram of type $5A_4$,
and hence it is isomorphic to the semi-direct product 
$\{\pm 1\}^5\semidirectproduct S_5$.
Note that $\Aut (F, h)$ acts on the dual lattice $(S_0)\dual$ 
of $S_0$ in a natural way,
and hence it  acts on the set of even overlattices of $S_0$.
We classify all  even overlattices of $S_0$
with the properties (ii) and (iii) in Claim~\ref{claim:overlattices}
up to the action of $\Aut (F, h)$.
The main tool is Nikulin's theory of discriminant forms 
of even lattices~\cite{MR525944}.
\par
\medskip
The set $F\cup \{h, l\}$ of vectors
form a basis of $S_0$.
Let
$$
(e_i^{(j)})\dual\quad(i=1, \dots, 4, \;j=1, \dots, 5),
\quad h\dual\quand  l\dual
$$
be the basis of  $(S_0)\dual$ dual to $F\cup \{h, l\}$.
We denote by $G$ the discriminant group $(S_0)\dual/S_0$ of $S_0$,
and by 
$$
\pr: (S_0)\dual\to G
$$
the natural projection.
Then $G$ is isomorphic to  $\F_5^{\oplus 5}\oplus \F_5$
with basis
$$
\pr((e_1^{(1)})\dual),
\;\;
\dots
\;\;, \;\;
\pr((e_1^{(5)})\dual),
\;\;
\pr(h\dual).
$$
With respect to this basis,
we denote the elements of $G$ by
$[x_1, \dots, x_5\,|\,y\,]$
with $x_1, \dots, x_5, y\in \F_5$.
The discriminant form $q: G\to\Q/2\Z$ of $S_0$ is given by
$$
q([x_1, \dots, x_5\,|\,y\,])=-\frac{4}{5}(x_1^2+\dots+x_5^2)+\frac{2}{5}y^2\; \bmod 2\Z
$$
The action of $\Aut (F, h)$ on $G=\F_5^{\oplus 5}\oplus \F_5$
is generated by the multiplications by $-1$ on $x_i$,
and the permutations of $x_1, \dots, x_5$.
We define subgroups $H_0, \dots, H_8$ of $G$ by their generators as follows:
\begin{eqnarray*}
H_0 &:=& \{0\}, \\
H_1 &:=& \ang{\,[0,0,2,2,2 \,|\, 2\,]\,}, \\
H_2 &:=& \ang{\,[2,2,2,2,2 \,|\, 0\,]\,}, \\
H_3 &:= &\ang{\,[0,1,2,2,2 \,|\,1\,]\,}, \\
H_4 &:= &\ang{\,[1,2,2,2,2 \,|\,2\,]\,},\\
H_5 &:=&\ang{\,[0,1,1,2,2\,|\,0\,]\,},\\
H_6 &:=& \ang{\,[1,0,1,2,2\,|\,0\,]\,, \,[0,1,2,1,3\,|\,0\,]\,},\\
H_7&:=&\ang{\,[1,0,0,1,1\,|\,1\,]\,,\,[0,1,1,1,3\,|\,3\,]\,},\\
H_8&:=&\ang{\,[1,0,1,1,2\,|\,2\,]\,,\,[0,1,1,3,3\,|\,0\,]\,}.
\end{eqnarray*}
We then put
$$
S_i:=\pr\inv(H_i)\;\;\subset\;\;(S_0)\dual.
$$
\begin{proposition}\label{prop:Hs} 
The submodules $S_0, \dots, S_8$ of $(S_0)\dual$
are even overlattices of $S_0$
with the properties {\rm (ii)} and {\rm (iii)} in Claim~\ref{claim:overlattices}.
The discriminant of $S_i$ is $-5^6$ for $i=0$, $-5^4$ for $i=1, \dots, 5$, and
$-5^2$ for $i=6, \dots, 8$.

Conversely,
if $S$ is an even overlattice 
of $S_0$
with the properties {\rm (ii)} and {\rm (iii)},
then there exists a unique $S_i$ among $S_0, \dots, S_8$ 
such that $S=g(S_i)$ holds for some $g\in \Aut (F,h)$.
\end{proposition}
\begin{proof}
The mapping $S\mapsto S/S_0$ gives rise to a one-to-one
correspondence between the set  of even overlattices $S$ of $S_0$ 
and  the set of
totally isotropic subgroups $H$ of $(G, q)$.
The inverse mapping is given by $H\mapsto \pr\inv (H)$.
If $\dim_{\F_5} H=d$, 
then the discriminant of $\pr\inv (H)$ is equal to $-5^{6-2d}$
(see Nikulin~\cite{MR525944}.)
\par
\medskip
For $v=[x_1, \dots, x_5\,|\,y\,]\in G$,
we put
$$
\delta(v):=(a,b,y)\in \Z_{\ge 0}\times \Z_{\ge 0}\times \F_5,
$$
where $a$ is the number of $\pm 1\in \F_5$ among $x_1, \dots, x_5$
and 
$b$ is the number of $\pm 2\in \F_5$ among $x_1, \dots, x_5$. 
Note that $\delta (v)=\delta (w)$ holds if and only if 
there exists $g\in \Aut (F,h)$ such that $g(v)=w$.
\begin{table}
\begin{tabular}{cccc}
the $(a,b,y)$-type  & the roots in $h\sperp$ & the set $E$ \\
\hline
$(0,0,0)$ & $5A_4$ & empty &$*$\\
$(0,2,\pm 1)$ & $A_9+3A_4$ & empty &\\
$(0,3,\pm 2)$ & $5A_4$ & empty&$*$ \\
$(0,5,0)$ & $5A_4$ & empty &$*$\\
$(1,1,0)$ & $E_8+3A_4$ & empty &\\
$(1,3,\pm 1)$ & $5A_4$ & empty&$*$ \\
$(1,4,\pm 2)$ & $5A_4$ & empty &$*$\\
$(2,0,\pm 2)$ & $A_9+3A_4$ & empty &\\
$(2,2,0)$ & $5A_4$ & empty &$*$\\
$(3,0,\pm 1)$ & $5A_4$ & empty &$*$\\
$(3,1,\pm 2)$ & $5A_4$ & empty &$*$\\
$(4,1,\pm 1)$ & $5A_4$ & empty &$*$\\
$(5,0,0)$ & $5A_4$ & empty &$*$\\
\end{tabular}
\vskip 3mm
\caption{The isotropic vectors in $(G, q)$}\label{table:isotr}
\end{table}
A vector $v\in G$ is isotropic with respect to $q$
if and only if $\delta(v)$ appears in the first column 
of Table~\ref{table:isotr}.
For each $(a,b,y)$-type  $\alpha$ in Table~\ref{table:isotr},
we choose a vector $v\in G$ such that $\delta(v)=\alpha$,
and calculate the even overlattice
$$
S_\alpha:=\pr\inv (\ang{v})
$$
of $S_0$.
The second column of Table~\ref{table:isotr}
presents the Dynkin type of the root system
$\shortset{r\in S_\alpha}{rh=0, r^2=-2}$,
and the third column presents the set
$E:=\shortset{e\in S_\alpha}{eh=1, e^2=0}$.
Hence
we see that
the following two conditions on a subgroup $H$ of $G$ are equivalent:
\begin{itemize}
\item[(I)]
The corresponding  submodule $\pr \inv (H)$ of $ (S_0)\dual$ is
an even overlattice 
of $S_0$
with the properties {\rm (ii)} and {\rm (iii)} in Claim~\ref{claim:overlattices}.
\item[(II)] For any $v\in H$,
$\delta (v)$ is an $(a, b, y)$--type with $*$ in Table~\ref{table:isotr}.
\end{itemize}
Using a computer,
we make the complete list 
of subgroups of $G$ that satisfy the condition (II) up to the action of $\Aut (F,h)$.
The complete set of representatives  is $\{H_0, \dots, H_8\}$ above.
\end{proof}
\begin{remark}
Since there  exist  no even unimodular lattices of signature $(1, 21)$
(see Serre~\cite[Theorem 5 in Chapter V]{MR0344216}),
all totally isotropic subgroups of $(G, q)$ are of dimension $\le 2$ over $\F_5$.
\end{remark}

\def\cprime{$'$} \def\cprime{$'$} \def\cprime{$'$} \def\cprime{$'$}
\providecommand{\bysame}{\leavevmode\hbox to3em{\hrulefill}\thinspace}
\providecommand{\MR}{\relax\ifhmode\unskip\space\fi MR }
\providecommand{\MRhref}[2]{%
  \href{http://www.ams.org/mathscinet-getitem?mr=#1}{#2}
}
\providecommand{\href}[2]{#2}

\end{document}